\documentclass[fleqn]{article}
\usepackage{graphicx} 
\usepackage{authblk}

\usepackage[margin=1in]{geometry}
\usepackage{xcolor}
\usepackage[hypertexnames=false,pdftex]{hyperref}

\definecolor{cadmiumgreen}{rgb}{0.0, 0.42, 0.24}
\definecolor{hotpink}{rgb}{1.0, 0.41, 0.71}
\definecolor{cgreen}{RGB}{0, 180, 100}
\hypersetup{
    colorlinks=true,
    linkcolor=cadmiumgreen,
    citecolor=cgreen,
    urlcolor=blue,
    linktocpage
}

\makeatletter
\let\newtitle\@title
\let\newauthor\@author
\let\newdate\@date
\makeatother
\usepackage{fancyhdr}
\makeatletter
\def\@maketitle{%
	\newpage
	\begin{center}%
		\let \footnote \thanks
		{\LARGE\bfseries \@title \par}%
		\vskip 2.5em%
			{\large
				\lineskip .5em%
				\begin{tabular}[t]{c}%
					\@author
				\end{tabular}\par}%
		\vskip 1em%
			{\large \@date}%
	\end{center}%
	\par
	\vskip 1.5em}
\makeatother

\AtBeginDocument{
	
	\renewcommand{\today}{}
}

\title{CBX: Python and Julia packages for\protect\\consensus-based interacting particle methods}

\author[1,*]{Rafael Bailo}
\affil[1]{Mathematical Institute, University of Oxford}
\author[2]{Alethea Barbaro}
\affil[2]{Technische Universiteit Delft}
\author[3]{Susana N. Gomes}
\affil[3]{Mathematics Institute, University of Warwick}
\author[4,5]{Konstantin Riedl}
\affil[4]{Technical University of Munich}
\affil[5]{Munich Center for Machine Learning}
\author[6,*]{\authorcr Tim Roith}
\affil[6]{Helmholtz Imaging, Deutsches Elektronen-Synchrotron DESY, Notkestr. 85, 22607 Hamburg, Germany}
\author[7]{Claudia Totzeck}
\affil[7]{University of Wuppertal}
\author[8,9]{Urbain Vaes}
\affil[8]{MATHERIALS team, Inria Paris}
\affil[9]{École des Ponts}

\affil[*]{Corresponding authors: bailo@maths.ox.ac.uk, tim.roith@desy.de}

\date{\today}

\usepackage{amsmath}
\usepackage{amssymb}
\usepackage{amsfonts}
\usepackage{hyperref}
\usepackage[maxbibnames=99]{biblatex}
\usepackage{listings}

\definecolor{codegreen}{rgb}{0,0.6,0}
\definecolor{codegray}{rgb}{0.5,0.5,0.5}
\definecolor{codepurple}{rgb}{0.58,0,0.82}
\definecolor{backcolour}{rgb}{0.95,0.95,0.92}

\makeatletter
\lst@CCPutMacro
    \lst@ProcessOther {"2A}{%
      \lst@ttfamily 
         {\raisebox{2pt}{*}}
         {\raisebox{1pt}{*}}
      }
    \@empty\z@\@empty
\makeatother

\lstdefinestyle{python}{
    commentstyle=\color{codegray},
    keywordstyle=\color{codegreen},
    numberstyle=\tiny\color{codegray},
    stringstyle=\color{codepurple},
    basicstyle=\ttfamily,
    breakatwhitespace=false,         
    breaklines=true,                 
    captionpos=b,                    
    keepspaces=true,
    showspaces=false,                
    showstringspaces=false,
    showtabs=false,                  
    tabsize=2
}

\begin{document}

\maketitle

\section{Summary}\label{summary}

We introduce \href{https://pdips.github.io/CBXpy/}{CBXPy} and
\href{https://pdips.github.io/ConsensusBasedX.jl/}{ConsensusBasedX.jl},
Python and Julia implementations of consensus-based interacting particle
systems (CBX), which generalise consensus-based optimization methods
(CBO) for global, derivative-free optimisation. The \emph{raison d'être}
of our libraries is twofold: on the one hand, to offer high-performance
implementations of CBX methods that the community can use directly,
while on the other, providing a general interface that can accommodate
and be extended to further variations of the CBX family. Python and
Julia were selected as the leading high-level languages in terms of
usage and performance, as well as for their popularity among the
scientific computing community. Both libraries have been developed with
a common \emph{ethos}, ensuring a similar API and core functionality,
while leveraging the strengths of each language and writing idiomatic
code.

\section{Mathematical background}\label{mathematical-background}

Consensus-based optimisation (CBO) is an approach to solve, for a given
(continuous) \emph{objective function}
\(f:\mathbb{R}^d \rightarrow \mathbb{R}\), the \emph{global minimisation
problem}

\[
x^* = \operatorname*{argmin}_{x\in\mathbb{R}^d} f(x),
\]

i.e., the task of finding the point \(x^*\) where \(f\) attains its
lowest value. Such problems arise in a variety of disciplines including
engineering, where \(x\) might represent a vector of design parameters
for a structure and \(f\) a function related to its cost and structural
integrity, or machine learning, where \(x\) could comprise the
parameters of a neural network and \(f\) the empirical loss, which
measures the discrepancy of the neural network prediction with the
observed data.

In some cases, so-called \emph{gradient-based methods} (those that
involve updating a guess of \(x^*\) by evaluating the gradient
\(\nabla f\)) achieve state-of-the-art performance in the global
minimisation problem. However, in scenarios where \(f\) is
\emph{non-convex} (when \(f\) could have many \emph{local minima}),
where \(f\) is \emph{non-smooth} (when \(\nabla f\) is not
well-defined), or where the evaluation of \(\nabla f\) is impractical
due to cost or complexity, \emph{derivative-free} methods are needed.
Numerous techniques exist for derivative-free optimisation, such as
\emph{random} or \emph{pattern search}
\autocite{friedman1947planning,rastrigin1963convergence,hooke1961direct},
\emph{Bayesian optimisation} \autocite{movckus1975bayesian} or
\emph{simulated annealing} \autocite{henderson2003theory}. Here, we
focus on \emph{particle-based methods}, specifically, consensus-based
optimisation (CBO), as proposed by \textcite{pinnau2017consensus}, and
the consensus-based taxonomy of related techniques, which we term
\emph{CBX}.

CBO uses a finite number \(N\) of \emph{agents} (particles),
\(x_t=(x_t^1,\dots,x_t^N)\), dependent on time \(t\), to explore the
landscape of \(f\) without evaluating any of its derivatives (as do
other CBX methods). The agents evaluate the objective function at their
current position, \(f(x_t^i)\), and define a \emph{consensus point}
\(c_\alpha\). This point is an approximation of the global minimiser
\(x^*\), and is constructed by weighing each agent's position against
the \emph{Gibbs-like distribution} \(\exp(-\alpha f)\)
\autocite{boltzmann1868studien}. More rigorously,

\[
c_\alpha(x_t) =
\frac{1}{ \sum_{i=1}^N \omega_\alpha(x_t^i) }
\sum_{i=1}^N x_t^i \, \omega_\alpha(x_t^i),
\quad\text{where}\quad
\omega_\alpha(\,\cdot\,) = \mathrm{exp}(-\alpha f(\,\cdot\,)),
\]

for some \(\alpha>0\). The exponential weights in the definition favour
those points \(x_t^i\) where \(f(x_t^i)\) is lowest, and comparatively
ignore the rest, particularly for larger \(\alpha\). If all the found
values of the objective function are approximately the same,
\(c_\alpha(x_t)\) is roughly an arithmetic mean. Instead, if one
particle is much better than the rest, \(c_\alpha(x_t)\) will be very
close to its position.

Once the consensus point is computed, the particles evolve in time
following the \emph{stochastic differential equation} (SDE)

\[
\mathrm{d}x_t^i =
-\lambda\ \underbrace{
\left( x_t^i - c_\alpha(x_t) \right) \mathrm{d}t
}_{
\text{consensus drift}
}
+ \sigma\ \underbrace{
\left\| x_t^i - c_\alpha(x_t) \right\| \mathrm{d}B_t^i
}_{
\text{scaled diffusion}
},
\]

where \(\lambda\) and \(\sigma\) are positive parameters, and where
\(B_t^i\) are independent Brownian motions in \(d\) dimensions. The
\emph{consensus drift} is a deterministic term that drives each agent
towards the consensus point, with rate \(\lambda\). Meanwhile, the
\emph{scaled diffusion} is a stochastic term that encourages exploration
of the landscape. The scaling factor of the diffusion is proportional to
the distance of the particle to the consensus point. Hence, whenever the
position of a particle and the location of the weighted mean coincide,
the particle stops moving. On the other hand, if the particle is far
away from the consensus, its evolution has a stronger exploratory
behaviour. While both the agents' positions and the consensus point
evolve in time, it has been proven that all agents eventually reach the
same position and that the consensus point \(c_\alpha(x_t)\) is a good
approximation of \(x^*\)
\autocite{carrillo2018analytical,fornasier2021consensus}. Other
variations of the method, such as CBO with anisotropic noise
\autocite{carrillo2021consensus}, \emph{polarised CBO}
\autocite{bungert2022polarized}, or \emph{consensus-based sampling}
(CBS) \autocite{carrillo2022consensus} have also been proposed.

In practice, the solution to the SDE above cannot be found exactly.
Instead, an \emph{Euler--Maruyama scheme} \autocite{KP1992} is used to
update the position of the agents. The update is given by

\[
x^i \gets x^i
-\lambda \,\Delta t
\left( x^i - c_\alpha(x) \right)
+ \sigma\sqrt{\Delta t}\
\left\| x^i - c_\alpha(x) \right\| \xi^i,
\]

where \(\Delta t > 0\) is the \emph{step size} and
\(\xi^i \sim \mathcal{N}(0,\mathrm{Id})\) are independent, identically
distributed, standard normal random vectors.

As a particle-based family of methods, CBX is conceptually related to
other optimisation approaches which take inspiration from biology, like
\emph{particle-swarm optimisation} (PSO) \autocite{kennedy1995particle},
from physics, like \emph{simulated annealing} (SA)
\autocite{henderson2003theory}, or from other heuristics
\autocite{mohan2012survey,karaboga2014comprehensive,yang2009firefly,bayraktar2013wind}.
However, unlike many such methods, CBX has been designed to be
compatible with rigorous convergence analysis at the mean-field level
\autocite[the infinite-particle limit, see][]{huang2021MFLCBO}. Many
convergence results have been shown, whether in the original formulation
\autocite{carrillo2018analytical,fornasier2021consensus}, for CBO with
anisotropic noise
\autocite{carrillo2021consensus,fornasier2021convergence}, with memory
effects \autocite{riedl2022leveraging}, with truncated noise
\autocite{fornasier2023consensus}, for polarised CBO
\autocite{bungert2022polarized}, and PSO
\autocite{qiu2022PSOconvergence}. The relation between CBO and
\emph{stochastic gradient descent} has been recently established by
\textcite{riedl2023gradient}, which suggests a previously unknown yet
fundamental connection between derivative-free and gradient-based
approaches.

\begin{figure}
\centering
\includegraphics[width=.7\textwidth]{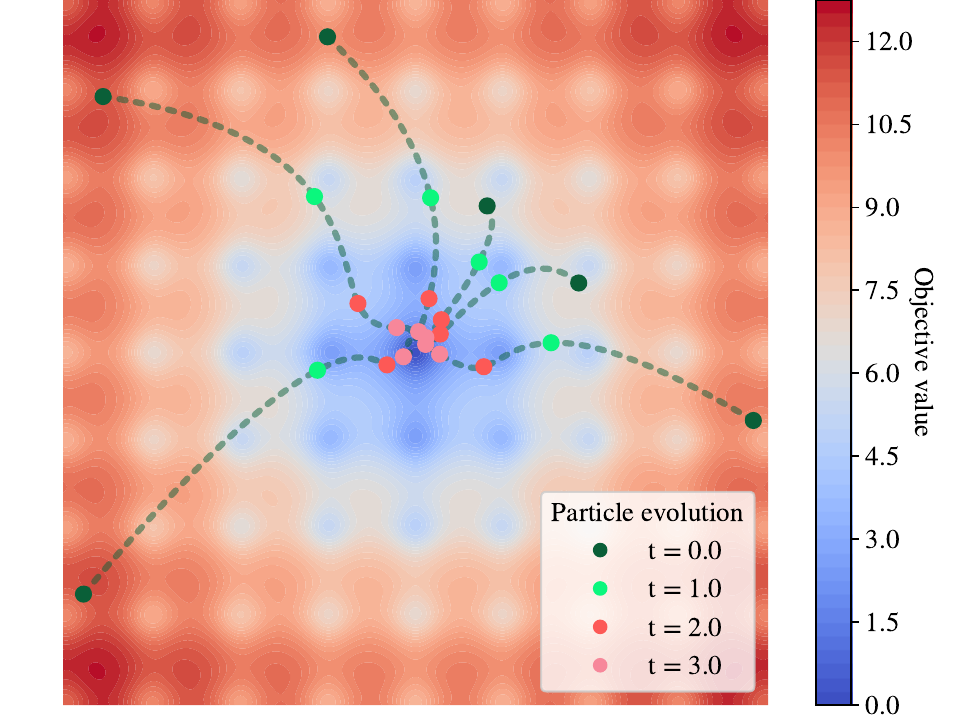}
\caption{Typical evolution of a CBO method minimising the Ackley
function \autocite{ackley2012connectionist}.}
\end{figure}

CBX methods have been successfully applied and extended to several
different settings, such as constrained optimisation problems
\autocite{fornasier2020consensus_sphere_convergence,borghi2021constrained},
multi-objective optimisation
\autocite{borghi2022adaptive,klamroth2022consensus}, saddle-point
problems \autocite{huang2022consensus}, federated learning tasks
\autocite{carrillo2023fedcbo}, uncertainty quantification
\autocite{althaus2023consensus}, or sampling
\autocite{carrillo2022consensus}.

\section{Statement of need}\label{statement-of-need}

In general, very few implementations of CBO already exist, and none have
been designed with the generality of other CBX methods in mind. Here, we
summarise the related software:

Regarding Python, we refer to \texttt{PyPop7} \autocite{duan2023pypop7}
and \texttt{scikit-opt} \autocite{scikitopt} for a collection of various
derivative-free optimisation strategies. For packages connected to
Bayesian optimisation, we refer to \texttt{BayesO} \autocite{Kim2023},
\texttt{bayesian-optimization} \autocite{Bayesian14}, \texttt{GPyOpt}
\autocite{gpyopt2016}, \texttt{GPflowOpt} \autocite{GPflowOpt2017},
\texttt{pyGPGO} \autocite{Jiménez2017}, \texttt{PyBADS}
\autocite{Singh2024} and \texttt{BoTorch}
\autocite{balandat2020botorch}. Furthermore, CMA-ES
\autocite{hansen1996adapting} was implemented in \texttt{pycma}
\autocite{hansen2019pycma}. To the best of our knowledge the connection
between consensus-based methods and evolution strategies is not fully
understood, and is therefore an interesting future direction. PSO and SA
implementations are already available in \texttt{PySwarms}
\autocite{miranda2018pyswarms}, \texttt{scikit-opt}
\autocite{scikitopt}, \texttt{DEAP} \autocite{deapJMLR2012} and
\texttt{pagmo} \autocite{pagmo2017}. They are widely used by the
community and provide a rich framework for the respective methods.
However, adjusting these implementations to CBO is not straightforward.
The first publicly available Python packages implementing CBX algorithms
were given by some of the authors together with collaborators.
\textcite{Igor_CBOinPython} implement standard CBO
\autocite{pinnau2017consensus}, and the package \texttt{PolarCBO}
\autocite{Roith_polarcbo} provides an implementation of polarised CBO
\autocite{bungert2022polarized}.
\href{https://pdips.github.io/CBXpy/}{CBXPy} is a significant extension
of the latter, which was tailored to the polarised variant. The code
architecture was generalised, which allowed the implementation of the
whole CBX family within a common framework.

\sloppy
Regarding Julia, PSO and SA methods are, among others, implemented in
\texttt{Optim.jl} \autocite{mogensen2018optim},
\texttt{Metaheuristics.jl} \autocite{mejia2022metaheuristics}, and
\texttt{Manopt.jl} \autocite{Bergmann2022}. PSO and SA are also included
in the meta-library \texttt{Optimization.jl} \autocite{DR2023}, as well
as Nelder--Mead, which is a direct search method. The latter is also
implemented in \texttt{Manopt.jl} \autocite{Bergmann2022}, which further
provides a manifold variant of CMA-ES \autocite{colutto2009cma}. One of
the authors gave the first specific Julia implementation of standard CBO
\texttt{Consensus.jl} \autocite{Bailo_consensus}. That package has now
been deprecated in favour of
\href{https://pdips.github.io/ConsensusBasedX.jl/}{ConsensusBasedX.jl},
which improves the performance of the CBO implementation with a
type-stable and allocation-free implementation. The package also adds a
CBS implementation, and overall presents a more general interface that
accomodates the wider CBX class of methods.

\section{Features}\label{features}

\href{https://pdips.github.io/CBXpy/}{CBXPy} and
\href{https://pdips.github.io/ConsensusBasedX.jl/}{ConsensusBasedX.jl}
provide a lightweight and high-level interface. An existing function can
be optimised with just one call. Method selection, parameters, different
approaches to particle initialisation, and termination criteria can be
specified directly through this interface, offering a flexible point of
entry for the casual user. Some of the methods provided are standard CBO
\autocite{pinnau2017consensus}, CBO with mini-batching
\autocite{carrillo2021consensus}, polarised CBO
\autocite{bungert2022polarized}, CBO with memory effects
\autocite{grassi2020particle,riedl2022leveraging}, and consensus-based
sampling (CBS) \autocite{carrillo2022consensus}. Parallelisation tools
are available.

A more proficient user will benefit from the fully documented interface,
which allows the specification of advanced options (e.g., debug output,
the noise model, or the numerical approach to the matrix square root of
the weighted ensemble covariance matrix). Both libraries offer
performance evaluation methods as well as visualisation tools.

Ultimately, a low-level interface (including documentation and full-code
examples) is provided. Both libraries have been designed to express
common abstractions in the CBX family while allowing customisation.
Users can easily implement new CBX methods or modify the behaviour of
the existing implementation by strategically overriding certain hooks.
The stepping of the methods can also be controlled manually.

\subsection{CBXPy}\label{cbxpy}

\begin{figure}[!ht]
\centering
\includegraphics[height=0.06\textheight, ]{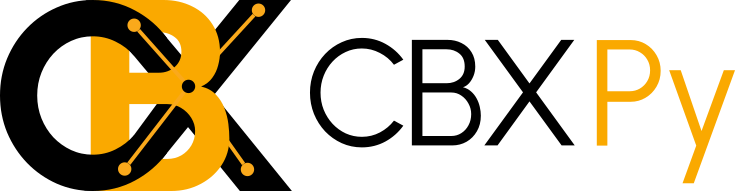}
\caption{CBXPy logo.}
\end{figure}

Most of the \href{https://pdips.github.io/CBXpy/}{CBXPy} implementation
uses basic Python functionality, and the agents are handled as an
array-like structure. For certain specific features, like
broadcasting-behaviour, array copying, and index selection, we fall back
to the \texttt{numpy} implementation \autocite{harris2020array}.
However, it should be noted that an adaptation to other array or tensor
libraries like PyTorch \autocite{paszke2019pytorch} is straightforward.
Compatibility with the latter enables gradient-free deep learning
directly on the GPU, as demonstrated in the documentation.\\

The library is available on
\href{https://github.com/pdips/CBXpy}{GitHub} and can be installed via
\texttt{pip}. It is licensed under the MIT license. Below, we provide a
short example on how to optimise a function with CBXPy.

\begin{lstlisting}[language=Python, style=python]
from cbx.dynamics import CBO        # import the CBO class
f = lambda x: x[0]**2 + x[1]**2     # define the function to minimise
x = CBO(f, d=2).optimize()          # run the optimisation
\end{lstlisting}

More examples and details on the implementation are available in the
\href{https://pdips.github.io/CBXpy/}{documentation}.

\subsection{ConsensusBasedX.jl}\label{consensusbasedx.jl}

\begin{figure}[!ht]
\centering
\includegraphics[height=0.06\textheight, ]{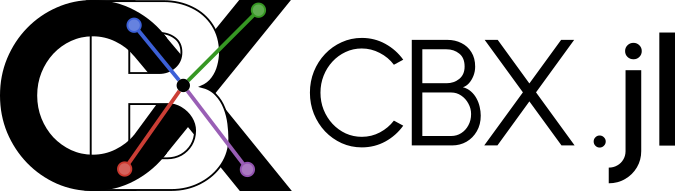}
\caption{ConsensusBasedX.jl logo.}
\end{figure}

\href{https://pdips.github.io/ConsensusBasedX.jl/}{ConsensusBasedX.jl}
has been almost entirely written in native Julia (with the exception of
a single call to LAPACK). The code has been developed with performance
in mind, thus the critical routines are fully type-stable and
allocation-free. A specific tool is provided to benchmark a typical
method iteration, which can be used to detect allocations. Through this
tool, unit tests are in place to ensure zero allocations in all the
provided methods. The benchmarking tool is also available to users, who
can use it to test their implementations of \(f\), as well as any new
CBX methods.

Basic function minimisation can be performed by running:

\begin{lstlisting}[language=Python, style=python]
using ConsensusBasedX               # load the ConsensusBasedX package
f(x) = x[1]^2 + x[2]^2              # define the function to minimise
x = minimise(f, D = 2)              # run the minimisation
\end{lstlisting}

The library is available on
\href{https://github.com/PdIPS/ConsensusBasedX.jl}{GitHub}. It has been
registered in the
\href{https://github.com/JuliaRegistries/General}{general Julia
registry}, and therefore it can be installed by running
\texttt{{]}add\ ConsensusBasedX}. It is licensed under the MIT license.
More examples and full instructions are available in the
\href{https://pdips.github.io/ConsensusBasedX.jl/}{documentation}.

\section{Acknowledgements}\label{acknowledgements}

We thank the Lorentz Center in Leiden for their kind hospitality during
the workshop ``Purpose-driven particle systems'' in Spring 2023, where
this work was initiated. RB was supported by the Advanced Grant
Nonlocal-CPD (Nonlocal PDEs for Complex Particle Dynamics: Phase
Transitions, Patterns and Synchronisation) of the European Research
Council Executive Agency (ERC) under the European Union's Horizon 2020
research and innovation programme (grant agreement No.~883363) and by
the EPSRC grant EP/T022132/1 ``Spectral element methods for fractional
differential equations, with applications in applied analysis and
medical imaging''. KR acknowledges support from the German Federal
Ministry of Education and Research and the Bavarian State Ministry for
Science and the Arts. TR acknowledges support from DESY (Hamburg,
Germany), a member of the Helmholtz Association HGF. This research was
supported in part through the Maxwell computational resources operated
at Deutsches Elektronen-Synchrotron DESY, Hamburg, Germany. UV
acknowledges support from the Agence Nationale de la Recherche under
grant ANR-23-CE40-0027 (IPSO).

\sloppy 
\printbibliography

\end{document}